\documentclass{amsart}[12pt]
\usepackage{amssymb}
\usepackage{amsfonts}
\usepackage{latexsym}

\newtheorem{theorem}{Theorem}[section]
\newtheorem{lemma}[theorem]{Lemma}
\newtheorem{proposition}[theorem]{Proposition}
 
\theoremstyle{definition}

\newtheorem{conjecture}[theorem]{Conjecture}

\newcommand{\g}{\mathfrak{g}}

\newcommand{\ben}{\begin{enumerate}}
\newcommand{\een}{\end{enumerate}}

\begin{document}

\title{On the Cachazo-Douglas-Seiberg-Witten conjecture for 
simple Lie algebras}

\begin{abstract} 
Recently, motivated by supersymmetric gauge theory, 
Cachazo, Douglas, Seiberg, and Witten proposed a conjecture 
about finite dimensional simple Lie
algebras, and checked it in the classical cases.
We prove the conjecture for type $G_2$, and also verify
a consequence of the conjecture in the general case. 
\end{abstract}

\author{Pavel Etingof}
\address{Department of Mathematics, Massachusetts Institute of Technology,
Cambridge, MA 02139, USA}
\email{etingof@math.mit.edu}

\author{Victor Kac}
\address{Department of Mathematics, Massachusetts Institute of Technology,
Cambridge, MA 02139, USA}
\email{kac@math.mit.edu}

\maketitle

\section{The CDSW conjecture}

Let $\g$ be a simple finite dimensional Lie algebra over ${\Bbb C}$.
We fix an invariant form on $\g$ and do not distinguish $\g$ and
$\g^*$. Let $g$ be the dual Coxeter number of $\g$.  
Consider the associative algebra $R=\wedge(\g\oplus \g)$. This algebra 
is naturally $\Bbb Z_+$-bigraded 
(with the two copies of $\g$ sitting in degrees (1,0) and (0,1), 
respectively). The degree (2,0),(1,1), and (0,2) components 
of $R$ are $\wedge^2\g,\g\otimes \g,\wedge^2\g$, 
respectively, hence each of them canonically contains a copy of $\g$. 
Let $I$ be the ideal in $R$ generated by these three copies of $\g$, and 
$A=R/I$. 

The associative algebra $A$ may be interpreted as follows.  
Let $\Pi \g$ denote $\g$ regarded as an odd vector space.
Then $R$ may be thought of as the algebra of regular functions on 
$\Pi\g\times \Pi\g$. 
We have the supercommutator map 
$\lbrace{,\rbrace}: \Pi\g\times \Pi\g\to \g$ 
given by the formula $\lbrace{X,Y\rbrace}=XY+YX$,
where the products are taken in the universal enveloping algebra
(this is a morphism of supermanifolds). The ideal $I$ in $R$ 
is then defined by the relations $\lbrace{X,X\rbrace}=0$, 
$\lbrace{X,Y\rbrace}=0$, $\lbrace{Y,Y\rbrace}=0$, so $A$ 
is the algebra of functions on the superscheme defined by these equations. 

In [CDSW,W], the following conjecture is proposed, and proved
for classical $\g$:

\begin{conjecture}\label{witten}
(i) The algebra $A^\g$ of $\g$-invariants in the algebra $A$ 
is generated by the unique invariant element $S$ of $A$ of degree (1,1)
(namely, $S=Tr|_V(XY)$, where $V$ is a non-trivial irreducible 
finite dimensional representation 
of $\g$).

(ii) $S^g=0$.

(iii) $S^{g-1}\ne 0$. 
\end{conjecture}

Here we prove the conjecture for $\g$
of type $G_2$. We believe that the method of the proof should be
relevant in the general case. 

We also prove in general the following result, which is a 
consequence of the conjecture. 

\begin{proposition}\label{triv}
Any homogeneous element of $A^\g$ 
is of degree $(d,d)$ for some $d$; 
therefore, the algebra 
$A^\g$ is purely even, and the natural action 
of $sl(2)$ on it (by linear transformations of $X,Y$) 
is trivial. 
\end{proposition}

{\bf Remark.} Conjecture \ref{witten} has the following
cohomological interpretation. Let $\g [x,y]$ denote 
the Lie algebra of polynomials of $x,y$ with values in $\g$. 
Coniser the algebra of relative cohomology 
$H^*(\g[x,y],\g,\Bbb C)$. It is graded by the cohomological degree 
$D$ and the $x,y$-degree $d$. It is clear that 
for any nonzero cochain, one has $D\le d$. Thus the homogeneous elements 
of the algebra $H^*(\g[x,y],\g,\Bbb C)$ for which $D=d$ form 
a subalgebra $E$. Conjecture \ref{witten} then states that 
$E$ is generated by an element $S$ of degree $2$ with the defining 
relation $S^g=0$. 

Note that in this formulation, one can obviously replace $\g [x,y]$ with its
quotient by any ideal spanned by homogeneous elements 
of $x,y$-degree $d\ge 3$. 

\section{Proof}
The main idea of the proof is to consider first the algebra 
$B$ of functions on the superscheme of $X\in \Pi \g$ such that 
$\lbrace{X,X\rbrace}=0$. This algebra was intensively studied by Kostant,  
Peterson, and others (see [K]) 
and is rather well understood. In particular, it is known 
that as a $\g$-module, $B$ is a direct sum of $2^{{\rm
rk}(\g)}$ non-isomorphic simple 
$\g$-modules $V_{\frak a}$ parametrized by abelian ideals
${\frak a}$ in a Borel subalgebra ${\frak b}$ of $\g$. Namely, if
${\frak a}$ 
is such an ideal of dimension $d$ then it defines (by taking its top exterior 
power) a nonzero vector $v_{\frak a}$ 
in $\wedge^d\g$ (defined up to scaling). This vector generates an irreducible
submodule in $\wedge^d\g$ (with highest weight vector $v_{\frak a}$), 
and the sum of these submodules is a complement to the kernel 
of the projection $\wedge \g\to B$. In fact, this is the unique 
invariant complement, because in each degree $d$ 
it is the eigenspace of the quadratic Casimir $C$ with eigenvalue $d$ (the
largest possible eigenvalue on $\wedge^d\g$; here the Casimir
is normalized so that $C|_{\g}=1$). Thus the (direct) sum 
of $V_{\frak a}$ is canonically identified with $B$ as a $\g$-module. 

The algebra $A$ is obtained from $B\otimes B$ by 
taking a quotient by the additional relation $\lbrace{X,Y\rbrace}=0$ 
and then taking invariants. In other words, $A=(B\otimes
B)^\g/L$, where 
$L$ is the space of invariants in the ideal $\tilde L$ in
$B\otimes B$ 
given by the relation $\lbrace{X,Y\rbrace}=0$. 

Thus, let us look more carefully at the algebra $(B\otimes
B)^\g$. From the above 
we see that a basis of $(B\otimes B)^\g$ is 
given by the elements $z_{\frak a}$, the canonical elements
in $V_{\frak a}\otimes V_{\frak a}^*$. The element $z_{\frak a}$ 
sits in bidegree $(d,d)$, where $d$ is the dimension of ${\frak
a}$. This implies Proposition \ref{triv} (the $sl(2)$ action on
$A^\g$ is trivial because the $sl(2)$-weight of a
vector in $A$ of degree $(d_1,d_2)$ is $d_1-d_2$).  

Now we proceed to prove the conjecture for $G_2$. 
Let $\omega_1$ and $\omega_2=\theta$ be the fundamental weights
of $G_2$ ($\theta$ is the highest root), and $\alpha_1,\alpha_2$ 
the corresponding simple roots (so that $\theta=3\alpha_1+2\alpha_2$). 
The abelian ideals in the corresponding Borel subalgebra are 
${\frak a}_0=0$, ${\frak a}_1=\Bbb C e_{\theta}$, ${\frak
a}_2=\Bbb C e_\theta\oplus \Bbb C e_{3\alpha_1+\alpha_2}$,
${\frak a}_3=\Bbb C e_\theta\oplus \Bbb C e_{3\alpha_1+\alpha_2}\oplus \Bbb C 
e_{2\alpha_1+\alpha_2}$.
Thus the $\g$-module $B$ has 4 irreducible
components: $V_0=\Bbb C$ in degree $0$, $V_1=\g$ in degree $1$, 
$V_2=V(3\omega_1)$ in degree 2, and $V_3=V(\omega_1+2\omega_2)$
in degree 3. 

Let $S$ be the invariant element of degree (1,1) in $\wedge \g\otimes
\wedge \g$. Clearly the projection of powers of $S$ onto $B\otimes
B$ is nonzero in degree 0,1,2,3 (this follows from the fact 
that for each degree $d$ one has a canonical decomposition 
$\wedge^d\g=V_{\frak a_d}\oplus {\rm Ker}(\wedge^d\g\to B[d])$,
where $B[d]$ is the degree $d$ component of $B$). 

The dual Coxeter number of $G_2$ is $4$.
Thus, our job is just to show that the ideal generated by the
relation $\lbrace{X,Y\rbrace}=0$ 
in $B\otimes B$ contains no nonzero $\g$-invariants. 

A $\g$-invariant of degree $(d+1,d+1)$ in this ideal is 
a linear combination of elements of 
the form $C_w:=\sum f_{ijk}(x_i\otimes x_j)w(x_k)$, where 
$x_i$ is an orthonormal basis of $\g$, 
$f_{ijk}$ are the structure constants of $\g$ in this basis, and 
$w: \g\to {\rm End}(V_{\frak a})=V_{\frak a}\otimes V_{\frak
a}^*$ 
is a homomorphism of
representations (here for brevity ${\frak a}:={\frak a}_d$). 
It is easy to check that 
the only homomorphisms $w$ 
relevant to our situation (i.e., for $d=1,2$) are 
given by the action of $\g$ on $V_{\frak a}$.
\footnote{
This has to do with the fact that the highest weight of $V_2$
does not involve $\omega_2$. 
The situation is exactly the same for type $B_2$,
except that it is not true there that all $w$ come from the
action of $\g$; this results in $S^3=0$ for $B_2$, while $S^3\ne
0$ for $G_2$.} 
So the result follows from the following lemma:

\begin{lemma} If $w$ is the action map then $C_w$ is zero in 
$B\otimes B$. 
\end{lemma}

\begin{proof}
Let us regard $C_w$ as an operator on $\wedge^{d+1}(\g)$. 
Then it has the form 
$$
C_w=\sum f_{ijk}W_{x_i}L_{x_k}I_{x_j},
$$
where $I,L,W$ are the operators of contraction, Lie derivative,
and wedging, respectively. This can be more shortly written as 
$$
C_w=\sum W_{x_i}L_{[x_i,x_j]}I_{x_j}.
$$
We claim that if $v\in V_{\frak a}$ then $C_wv$ is zero
(this would imply that $C_w$ is zero in $B\otimes B$). 
Set $y=p_1\wedge...\wedge p_{d+1}$. Then we get 
$$
C_wy=\sum_{i,j,r,s}(-1)^{r+s} x_i\wedge
([[x_i,p_r],p_s]-[[x_i,p_s],p_r])
\wedge p_1\wedge ...\hat p_r...\hat p_s...\wedge p_{d+1}=
$$
$$
\sum_{i,j,r,s}(-1)^{r+s} x_i\wedge
[[p_s,p_r]x_i]
\wedge p_1\wedge ...\hat p_r...\hat p_s...\wedge p_{d+1}.
$$
This implies the result, since $V_{\frak a}$ coincides with the
kernel of the map $\wedge^{d+1}\g\to \g\otimes \wedge^{d-1}\g$,
given by 
$$
y\mapsto 
\sum_{r,s}(-1)^{r+s} [p_s,p_r]\otimes 
p_1\wedge ...\hat p_r...\hat p_s...\wedge p_{d+1}.
$$
\end{proof}

\vskip .1in

\centerline{\bf References}
\vskip .05in

[CDSV]
Freddy Cachazo, Michael R. Douglas, Nathan Seiberg, Edward Witten,
Chiral Rings and Anomalies in Supersymmetric Gauge Theory,
hep-th/0211170.

[K]
Bertram Kostant,  
On $\bigwedge{\frak g}$ for a semisimple Lie algebra ${\frak g}$, as
an equivariant 
module over the symmetric algebra $S({\frak g})$.
Analysis on homogeneous spaces and representation theory of Lie
groups, 
Okayama--Kyoto (1997), 129--144, Adv. Stud. Pure Math., 26, 
Math. Soc. Japan, Tokyo, 2000.

[W] Edward Witten, 
Chiral Ring Of Sp(N) And SO(N) Supersymmetric Gauge Theory In Four Dimensions
hep-th/0302194.

\end{document}